\documentclass{amsart}
\usepackage{amsfonts,amsmath,amsthm,amssymb}

\newtheorem{theorem}{Theorem}

\newtheorem{corollary}{Corollary}
\newtheorem{proposition}{Proposition}

\thanks{The research was supported by OTKA  No.K68383 and FAPESP 07/54571-8}


\begin{document}
\title[Integral group ring of the Suzuki sporadic simple group]{Integral group ring of the Suzuki\\sporadic simple group}
\date{}

\author{V.A.~Bovdi, A.B.~Konovalov, E.N.~Marcos}

\address{V.A.~Bovdi
\newline Institute of Mathematics, University of Debrecen,
\newline P.O.  Box 12, H-4010 Debrecen, Hungary
\newline Institute of Mathematics and Informatics, College of Ny\'\i regyh\'aza,
\newline S\'ost\'oi \'ut 31/b, H-4410 Ny\'\i regyh\'aza, Hungary}
\email{vbovdi@math.klte.hu}

\address{A.B.~Konovalov
\newline School of Computer Science, University of St Andrews,
\newline Jack Cole Building, North Haugh, St Andrews, Fife, KY16 9SX, Scotland}
\email{konovalov@member.ams.org}

\address{E.N.~Marcos
\newline Departamento de Matem\'atica, Universidade de Sao Paulo
\newline Caixa Postal 66281 05315-970 - Sao Paulo, Brazil}
\email{enmarcos@ime.usp.br}

\begin{abstract}
Using the Luthar--Passi method, we investigate the classical
Zassenhaus conjecture for the normalized unit group of the
integral group ring of the Suzuki sporadic simple group
$\verb"Suz"$. As a consequence, for this group we confirm
the Kimmerle's conjecture on prime graphs.
\end{abstract}

\maketitle

%#########################################

\section{Introduction and main results}

Let $V(\mathbb Z G)$ be  the normalized unit group of the
integral group ring $\mathbb Z G$ of  a finite group $G$. A long-standing
conjecture  of H.~Zassenhaus {\bf (ZC)} says that every torsion unit
$u\in V(\mathbb ZG)$ is conjugate within the rational group algebra
$\mathbb Q G$ to an element in $G$ (see \cite{Zassenhaus}).

For finite simple groups the main tool for the investigation of the
Zassenhaus conjecture is the Luthar--Passi method, introduced in
\cite{Luthar-Passi} to solve it for $A_{5}$ and then applied in
\cite{Luthar-Trama} for the case of $S_{5}$.
Later M.~Hertweck extended
the Luthar--Passi method and applied it for the investigation
of the Zassenhaus conjecture for $PSL(2,p)$, $p=7,11,13$ in \cite{Hertweck1},
and for $PSL(2,9) \cong A_6$ in \cite{Hertweck5}.
The Luthar--Passi method proved to be useful for groups containing
non-trivial normal subgroups as well. For some recent results we refer to
\cite{Bovdi-Hofert-Kimmerle,Bovdi-Konovalov, Hertweck2, Hertweck1,
Hertweck3, Hofert-Kimmerle}.
Also, some related properties and  some weakened variations of the
Zassenhaus  conjecture can be found in
\cite{Artamonov-Bovdi,Bleher-Kimmerle,Luthar-Trama}.

First of all, we need to introduce some notation. By $\# (G)$ we
denote the set of all primes dividing the order of $G$. The
Gruenberg--Kegel graph (or the prime graph) of $G$ is the graph
$\pi (G)$ with vertices labeled by the primes in $\# (G)$ and with
an edge from $p$ to $q$ if there is an element of order $pq$ in
the group $G$. In \cite{Kimmerle} W.~Kimmerle   proposed the
following weakened variation of the Zassenhaus conjecture:

\begin{itemize}
\item[]{\bf (KC)} \qquad
If $G$ is a finite group then $\pi (G) =\pi (V(\mathbb Z G))$.
\end{itemize}

In particular, in the same  paper  W.~Kimmerle verified   that
{\bf (KC)} holds for finite Frobenius and solvable groups. We remark
that with respect to the so-called $p$-version of the Zassenhaus
conjecture the investigation
of Frobenius groups was completed by  M.~Hertweck and the first
author  in \cite{Bovdi-Hertweck}. In \cite{Bovdi-Jespers-Konovalov,
Bovdi-Konovalov, Bovdi-Konovalov-McL, Bovdi-Konovalov-Ru,
Bovdi-Konovalov-M24, Bovdi-Konovalov-HS,
Bovdi-Konovalov-M23, Bovdi-Konovalov-Linton,
Bovdi-Konovalov-Siciliano}
{\bf (KC)} was confirmed for the Mathieu simple groups $M_{11}$,
$M_{12}$, $M_{22}$, $M_{23}$, $M_{24}$, sporadic Janko simple groups
$J_1$, $J_{2}$ and $J_3$,
the Higman-Sims sporadic group $HS$,
the McLaughlin sporadic group $McL$ and
the Rudvalis sporadic group $Ru$.

Here we continue these investigations for the Suzuki simple group
$\verb"Suz"$. Although using the Luthar--Passi method we cannot prove the
rational conjugacy for torsion units of its integral group ring, our
main result gives a lot of information on partial augmentations of
these units. In particular, we confirm the Kimmerle's conjecture for
this group.

Let $G$ be the Suzuki sporadic simple group $\verb"Suz"$ and let
\[
\begin{split}
\mathcal{C} =\{
C_{1}, C_{2a}, \; C_{2b}, \; C_{3a}, \;&
C_{3b}, \; C_{3c}, \; C_{4a}, \; C_{4b}, \; C_{4c}, \;
C_{4d}, \; C_{5a}, \; C_{5b}, \; C_{6a},\\  \; C_{6b},
\; C_{6c}, \quad   &C_{6d}, \; C_{6e}, \; C_{7a}, \;
C_{8a}, \; C_{8b}, \; C_{8c}, \; C_{9a}, \; C_{9b}, \;
C_{10a}, \;\\  C_{10b}, \; C_{11a}, \;& C_{12a}, \;
C_{12b}, \; C_{12c}, \; C_{12d}, \;
C_{12e}, \; C_{13a}, \; C_{13b}, \;C_{14a}, \; \\
C_{15a}, \;& C_{15b}, \; C_{15c}, \; C_{18a}, \;
C_{18b}, \; C_{20a}, \; C_{21a}, \; C_{21b}, \; C_{24a}
\}
\end{split}
\]
be the collection of all its conjugacy classes, where the first
index denotes the order of the elements of this conjugacy class
and $C_{1}=\{ 1\}$. Suppose $u=\sum \alpha_g g \in V(\mathbb Z G)$
has finite order $k>1$. Denote by
$\nu_{nt}=\nu_{nt}(u)=\varepsilon_{C_{nt}}(u)=\sum_{g\in C_{nt}}
\alpha_{g}$ the partial augmentation of $u$ with respect to
$C_{nt}$.
From the Berman--Higman Theorem
(see \cite{Berman} and \cite{Sandling}, Ch.5, p.102)
one knows that
$\nu_1 =\alpha_{1}=0$, so
\begin{equation}\label{E:1}
\sum_{C_{nt}\in \mathcal{C}\setminus C_1} \nu_{nt}=1.
\end{equation}
Hence, for any character $\chi$ of $G$, we get that $\chi(u)=\sum
\nu_{nt}\chi(h_{nt})$, where $h_{nt}$ is a representative of the
conjugacy class $ C_{nt}$.

Our main result is the following theorem.

\begin{theorem}\label{T:1}
Let $G$ denote the Suzuki  sporadic simple group $\verb"Suz"$. Let
$u$ be a torsion unit of $V(\mathbb ZG)$ of order $|u|$ and denote by
\[
\begin{split}
\mathfrak{P}(u)=( \nu_{2a}, \; \nu_{2b}, \; \nu_{3a}, \;&
\nu_{3b}, \; \nu_{3c}, \; \nu_{4a}, \; \nu_{4b}, \; \nu_{4c}, \;
\nu_{4d}, \; \nu_{5a}, \; \nu_{5b}, \; \nu_{6a},\\  \; \nu_{6b},
\; \nu_{6c}, \quad   &\nu_{6d}, \; \nu_{6e}, \; \nu_{7a}, \;
\nu_{8a}, \; \nu_{8b}, \; \nu_{8c}, \; \nu_{9a}, \; \nu_{9b}, \;
\nu_{10a}, \;\\  \nu_{10b}, \; \nu_{11a}, \;& \nu_{12a}, \;
\nu_{12b}, \; \nu_{12c}, \; \nu_{12d}, \;
\nu_{12e}, \; \nu_{13a}, \; \nu_{13b}, \;\nu_{14a}, \; \\
\nu_{15a}, \;& \nu_{15b}, \; \nu_{15c}, \; \nu_{18a}, \;
\nu_{18b}, \; \nu_{20a}, \; \nu_{21a}, \; \nu_{21b}, \; \nu_{24a}
\;) \in \mathbb Z^{42}
\end{split}
\]
the tuple of partial augmentations of $u$. The following
properties hold.

\begin{itemize}

\item[(i)] There is no elements of orders  $22$, $26$, $33$, $35$,
$39$, $55$, $65$, $77$, $91$ and  $143$  in $V(\mathbb ZG)$.
Equivalently, if $|u|$ not in
\[
\begin{split}
\qquad \qquad
\{ \; 28, \; 30, \; & 36, \; 40, \; 42, \; 45, \; 56, \; 60, \; 63, \; 72, \\
& \; 84, \; 90, \; 120, \; 126, \; 168, \; 180, \; 252, \; 360, \; 504 \; \},
\end{split}
\]
then $|u|$ coincides with the order of some element $g \in G$.

\item[(ii)] If $|u| \in \{7, 11\}$ then $u$ is rationally
conjugate to some $g \in G$.

\item[(iii)] If $|u|=2$, then one has
\[
\begin{split}
( \nu_{2a}, \nu_{2b} ) \in \{ (4, -3 ), &( 3, -2 ), ( 2, -1 ), (
1,0 ),\\ &( 0, 1 ), ( -1, 2 ),  ( -2, 3 ), ( -3, 4 ) \; \}.
\end{split}
\]

\item[(iv)] If $|u|=5$, then one has
\[
\begin{split}
( \nu_{5a}, \nu_{5b} ) \in \{ ( 5, -4 ), &( 4, -3 ), ( 3, -2 ), (
2, -1 ),\\ &( 1, 0 ), ( 0, 1 ),  ( -1, 2 ), ( -2, 3 ), ( -3, 4 )\;
\}.
\end{split}
\]

\item[(v)]  If $|u|=13$, then one has
\[
\begin{split}
\{ \; ( \nu_{13a}, \nu_{13b} ) \; \mid \; -8 \le \nu_{13a} \le 9, \quad \nu_{13a}+\nu_{13b}=1 \; \}.
%( \nu_{13a}, \nu_{13b} ) \in \{( 9, -8 ),( 8, -7 ), ( 7, -6 ), (6,
%-5 ), ( 5, -4 ), ( 4, -3 ), ( 3, -2 ),&\\   ( 2, -1 ), ( 1, 0 ),
%(0, 1 ), ( -1, 2 ),  ( -2, 3 ),( -3, 4 ), ( -4, 5 ),&\\ ( -5, 6 ),
%(-6, 7 ), ( -7, 8 ), ( -8, 9 )&\;  \}.
\end{split}
\]

\end{itemize}
\end{theorem}

For the case of torsion units of order 3, using our implementation of the
Luthar--Passi method, which we intend to make available in the GAP package
LAGUNA \cite{LAGUNA}, we are able to compute the set of 104 tuples containing
(likely as a proper subset) possible tuples of partial augmentations,
listed in the Appendix.

As an immediate consequence of  part (i) of the Theorem we obtain

\begin{corollary} If $G = \verb"Suz"$ then
$\pi(G)=\pi(V(\mathbb ZG))$.
\end{corollary}

\section{Preliminaries}
The following result is a reformulation of
the Zassenhaus conjecture in terms of vanishing
of partial augmentations of torsion units.

\begin{proposition}\label{P:1}
(see \cite{Luthar-Passi} and
Theorem 2.5 in \cite{Marciniak-Ritter-Sehgal-Weiss})
Let $u\in V(\mathbb Z G)$
be of order $k$. Then $u$ is conjugate in $\mathbb
QG$ to an element $g \in G$ if and only if for
each $d$ dividing $k$ there is precisely one
conjugacy class $C$ with partial augmentation
$\varepsilon_{C}(u^d) \neq 0 $.
\end{proposition}

The next result now yield that several partial augmentations
are zero.

\begin{proposition}\label{P:2}
(see \cite{Hertweck2}, Proposition 3.1;
\cite{Hertweck1}, Proposition 2.2)
Let $G$ be a finite
group and let $u$ be a torsion unit in $V(\mathbb
ZG)$. If $x$ is an element of $G$ whose $p$-part,
for some prime $p$, has order strictly greater
than the order of the $p$-part of $u$, then
$\varepsilon_x(u)=0$.
\end{proposition}

The key restriction on partial augmentations is given
by the following result that is the cornerstone of
the Luthar--Passi method.

\begin{proposition}\label{P:3}
(see \cite{Hertweck1,Luthar-Passi}) Let either $p=0$ or $p$ a prime
divisor of $|G|$. Suppose
that $u\in V( \mathbb Z G) $ has finite order $k$ and assume $k$ and
$p$ are coprime in case $p\neq 0$. If $z$ is a complex primitive $k$-th root
of unity and $\chi$ is either a classical character or a $p$-Brauer
character of $G$, then for every integer $l$ the number
$$
\mu_l(u,\chi, p ) =
\textstyle\frac{1}{k} \sum_{d|k}Tr_{ \mathbb Q (z^d)/ \mathbb Q }
\{\chi(u^d)z^{-dl}\}
$$
is a non-negative integer.
\end{proposition}

Note that if $p=0$, we will use the notation $\mu_l(u,\chi
, * )$ for $\mu_l(u,\chi , 0)$.

Finally, we will use the well-known bound for
orders of torsion units.

\begin{proposition}\label{P:4}  (see  \cite{Cohn-Livingstone})
The order of a torsion element $u\in V(\mathbb ZG)$
is a divisor of the exponent of $G$.
\end{proposition}

\section{Proof of the Theorem}

%#########################################

In this section we denote by $G$ the Suzuki sporadic simple  group
$\verb"Suz"$.
The character table of $G$,
as well as the $p$-Brauer character tables, which will be denoted by
$\mathfrak{BCT}{(p)}$ where $p\in\{2,3,5,7,11,13\}$, can be found using
the computational algebra system GAP \cite{GAP},
which derives these data from \cite{AFG,ABC}. For the characters
and conjugacy classes we will use throughout the paper the same
notation, indexation inclusive, as used in the GAP Character Table Library.

It is well known (see \cite{GAP,Gorenstein}) that \quad $
|G|= 448345497600 = 2^{13} \cdot 3^{7} \cdot 5^{2} \cdot 7 \cdot
11 \cdot 13 $ and $exp(G)= 360360 = 2^{3} \cdot 3^{2} \cdot 5
\cdot 7 \cdot 11 \cdot 13$.

Since the group $G$ possesses elements of orders
2, 3, 4, 5, 6, 7, 8, 9, 10, 11, 12, 13, 14, 15, 18, 20, 21 and 24,
first of all we investigate units of some of these orders
(except for the units of orders 4, 6, 8, 9, 10, 12, 14, 15, 18, 20, 21 and 24).
After this, by Proposition \ref{P:4}, the order of each torsion unit divides
the exponent of $G$, so to prove the Kimmerle's conjecture, it remains to
consider units of orders 22, 26, 33, 35, 39, 55, 65, 77, 91 and 143.
We will prove that no units of all these orders do appear in $V(\mathbb ZG)$.
Since we omit orders  28, 30, 36, 40, 42, 45 and 63 that do not
contribute to {\bf(KC)}, we need to add to the list of exceptions in part
(i) of Theorem also orders 56, 60, 72, 84, 90, 120, 126, 168, 180,
252, 360 and 504, but no more because of restrictions imposed
by the exponent of G.

Now we consider separately each possible value of $|u|$.

\noindent $\bullet $ Let $|u|\in \{7,11\}$. Since there is only one conjugacy
class in $G$ consisting of elements or order $|u|$, this case follows at
once from Proposition \ref{P:2}. Thus, for units of orders $5$ and $7$ we
obtained that there is precisely one conjugacy class with non-zero partial
augmentation. Proposition \ref{P:1} then yields part (ii) of the Theorem.

\noindent $\bullet $ Let $u$ be an involution. By (\ref{E:1}) and
Proposition \ref{P:2} we get $ \nu_{2a}+\nu_{2b}=1$. Put $t_1=15
\nu_{2a} -  \nu_{2b}$ and $t_2=7\nu_{2a} -3\nu_{2b}$. Applying
Proposition \ref{P:3} to the ordinary character $\chi _{2}$ and
3-Brauer character $\chi_{3}$ we get the following system of inequalities
\[
\begin{split}
\mu_{0}(u,\chi_{2},*) & = \textstyle \frac{1}{2} (t_1+ 143) \geq
0; \quad \mu_{1}(u,\chi_{2},*) = \textstyle \frac{1}{2} (-t_1 + 143) \geq 0; \\ % 2
\mu_{0}(u,\chi_{3},3) & = \textstyle \frac{1}{2} (2t_2 + 78) \geq
0; \quad
\mu_{1}(u,\chi_{3},3)  = \textstyle \frac{1}{2} (-2t_2  + 78) \geq 0. \\ % 76
\end{split}
\]
From the requirement that all $\mu _{i}(u,\chi _{j},p)$ must be
non-negative integers it can be deduced that $(\nu_{2a},\nu_{2b})$
satisfies the conditions of part (iii) of the Theorem.

\noindent $\bullet $ Let $u$ be a unit of order $3$. By (\ref{E:1}) and
Proposition \ref{P:2} we obtain that
\begin{equation*}
\nu _{3a}+\nu _{3b}+\nu _{3c}=1.
\end{equation*}
Put $t_1 = 35 \nu_{3a} + 8 \nu_{3b} -  \nu_{3c}$ and $t_2 = 14
\nu_{3a} - 13 \nu_{3b} - 4 \nu_{3c}$. Again applying Proposition \ref{P:3}
to the ordinary characters $\chi_{2}$, $\chi_{3}$ and 2-Brauer characters
$\chi_{2}$ and $\chi_{7}$, we obtain the following system
of inequalities
\[
\begin{split}
\mu_{0}(u,\chi_{2},*) & = \textstyle \frac{1}{3} (2t_1 + 143) \geq
0; \quad \mu_{1}(u,\chi_{2},*)= \textstyle \frac{1}{3} (-t_1+143) \geq 0; \\ % 2
\mu_{0}(u,\chi_{3},*) & = \textstyle \frac{1}{3} (-t_2 + 364) \geq 0;
\quad \mu_{1}(u,\chi_{3},*)  = \textstyle \frac{1}{3} (t_2 + 364) \geq 0; \\ % 4
&\mu_{0}(u,\chi_{2},2) = \textstyle \frac{1}{3} (-50 \nu_{3a} + 4\nu_{3b} + 4 \nu_{3c} + 110) \geq 0; \\ % 71
&\mu_{0}(u,\chi_{7},2) = \textstyle \frac{1}{3} (142 \nu_{3a} -20 \nu_{3b} + 16 \nu_{3c} + 638) \geq 0, \\ % 77
\end{split}
\]
that has only $101$ non-trivial and three trivial solutions
$(\nu_{3a},\nu _{3b},\nu _{3c})$ which are listed in the Appendix.

\noindent$\bullet$ Let $u$ be a unit of order $5$. By (\ref{E:1})
and Proposition \ref{P:2} we have $ \nu_{5a}+\nu_{5b}=1$. Put $t_1
= 8 \nu_{5a} + 3 \nu_{5b}$. From the ordinary character table and
Brauer character tables for $p=2,3$ we obtain the following
system of inequalities
\[
\begin{split}
\mu_{0}(u,\chi_{2},*) & = \textstyle \frac{1}{5} (4 t_1+ 143) \geq
0; \quad \mu_{1}(u,\chi_{2},*) = \textstyle \frac{1}{5} (-t_1+ 143) \geq 0; \\ % 2
&\qquad  \mu_{0}(u,\chi_{2},2)  = \textstyle \frac{1}{5} (-20 \nu_{5a} + 110) \geq 0; \\ % 51
& \qquad \mu_{0}(u,\chi_{2},3)  = \textstyle \frac{1}{5} (16 \nu_{5a} - 4 \nu_{5b} + 64) \geq 0. \\ % 65
\end{split}
\]
From the requirement that all $\mu _{i}(u,\chi _{j},p)$ must be
non-negative integers it can be deduced that $(\nu_{5a},\nu_{5b})$
satisfies the condition of part (iv) of the Theorem.

\noindent$\bullet$ Let $u$ be a unit of order $13$. By (\ref{E:1})
and Proposition \ref{P:2} we have $\nu_{13a}+\nu_{13b}=1$. Put
$t_1 = 7 \nu_{13a} - 6 \nu_{13b}$. Applying Proposition \ref{P:3}
to the Brauer character tables for $p=2,3$, we get the following
system of inequalities
\[
\begin{split}
\mu_{1}(u,\chi_{2},2) & = \textstyle \frac{1}{13} (t_1 + 110) \geq
0; \quad \mu_{1}(u,\chi_{10},3)  = \textstyle \frac{1}{13} (-t_1 + 5103) \geq 0; \\ % 7
& \mu_{2}(u,\chi_{2},2) = \textstyle \frac{1}{13} (-6 \nu_{13a} + 7 \nu_{13b} + 110) \geq 0, \\ % 4
\end{split}
\]
which has integral solution $(\nu_{13a},\nu_{13b})$ listed in part
(v) of the Theorem.

\noindent$\bullet$ Let $u$ be a unit of order $22$. By (\ref{E:1})
and Proposition \ref{P:2} we have
$$
\nu_{2a}+\nu_{2b}+\nu_{11a}=1.
$$
Put $t_1 = 15\nu_{2a}-\nu_{2b}$ and $t_2 = 44 \nu_{2a}+\nu_{11a}$.
Since $|u^{11}|=2$, by part (iii) of the Theorem we have
eight cases, which we consider separately.

Case 1. Let $\chi(u^{11}) = \chi(2a)$. Using Proposition \ref{P:3}
for the characters $\chi_{2}, \chi_{3}$ and $\chi_{4}$ of
$G$, we get the following system of inequalities
\[
\begin{split}
\mu_{0}(u,\chi_{2},*) & = \textstyle \frac{1}{22} (10t_1+ 158)\geq 0; \quad \mu_{11}(u,\chi_{2},*) = \textstyle \frac{1}{22} (-10 t_1 + 128) \geq 0; \\ % 4
\mu_{0}(u,\chi_{3},*) & = \textstyle \frac{1}{22} (10 t_2+ 418)\geq 0; \quad \mu_{11}(u,\chi_{3},*)= \textstyle \frac{1}{22} (-10 t_2+ 330) \geq 0; \\ % 8
& \mu_{11}(u,\chi_{4},*)  = \textstyle \frac{1}{22} (-120 \nu_{2a} - 200 \nu_{2b} + 10 \nu_{11a} + 758) \geq 0; \\ % 12
&\quad \qquad \mu_{1}(u,\chi_{2},*) = \textstyle \frac{1}{22} (15 \nu_{2a} -  \nu_{2b} + 128) \geq 0; \\ % 2
&\quad \qquad \mu_{1}(u,\chi_{3},*) = \textstyle \frac{1}{22} (44 \nu_{2a} +  \nu_{11a} + 319) \geq 0. \\ % 6
\end{split}
\]

Case 2. Let $\chi(u^{11}) =\chi(2b)$. Using Proposition \ref{P:3}
for the characters $\chi_{2}, \chi_{3}$  and $\chi_{4}$ of
$G$, we obtain the following system
\[
\begin{split}
\mu_{0}(u,\chi_{2},*) & = \textstyle \frac{1}{22} (10t_1 + 142)\geq 0; \quad \mu_{11}(u,\chi_{2},*)= \textstyle \frac{1}{22} (-10t_1 + 144) \geq 0; \\ % 2
\mu_{1}(u,\chi_{2},*) & = \textstyle \frac{1}{22} (t_1+ 144) \geq 0; \qquad \mu_{1}(u,\chi_{3},*) = \textstyle \frac{1}{22} (t_2+ 363) \geq 0; \\ % 6
\mu_{0}(u,\chi_{3},*) & = \textstyle \frac{1}{22} (10 t_2+ 374)\geq 0; \quad \mu_{11}(u,\chi_{3},*)= \textstyle \frac{1}{22} (-10t_2+ 374) \geq 0; \\ % 8
&\mu_{0}(u,\chi_{4},*) = \textstyle \frac{1}{22} (10(12\nu_{2a} + 20 \nu_{2b} - \nu_{11a}) + 790) \geq 0; \\ % 9
&\mu_{11}(u,\chi_{4},*) = \textstyle \frac{1}{22} (-10(12 \nu_{2a} +20 \nu_{2b} -\nu_{11a}) + 750) \geq 0. \\ % 12
\end{split}
\]

Case 3. Let $\chi(u^{11}) = 4\chi(2a)-3\chi(2b)$.
Using Proposition \ref{P:3} for the characters $\chi_{2}$ and
$\chi_{3}$ of $G$, we obtain the following system
\[
\begin{split}
\mu_{0}(u,\chi_{2},*) & = \textstyle \frac{1}{22} (10 t_1+ 206)\geq 0; \quad \mu_{1}(u,\chi_{2},*) = \textstyle \frac{1}{22} (t_1+ 80) \geq 0; \\ % 2
\mu_{11}(u,\chi_{2},*) & = \textstyle \frac{1}{22} (-10 t_1+ 80)\geq 0; \quad \mu_{0}(u,\chi_{3},*) = \textstyle \frac{1}{22} (10 t_2+ 550) \geq 0; \\ % 5
\mu_{1}(u,\chi_{3},*) & = \textstyle \frac{1}{22} (t_2+ 187)
\geq0; \qquad\mu_{11}(u,\chi_{3},*) = \textstyle \frac{1}{22} (-10 t_2+ 198) \geq 0. \\ % 8
\end{split}
\]

Case 4. Let $\chi(u^{11}) = 3\chi(2a)-2\chi(2b)$. Using
Proposition \ref{P:3} for the characters $\chi_{2}$ and
$\chi_{3}$ of $G$, we obtain  the following system
\[
\begin{split}
\mu_{0}(u,\chi_{2},*) & = \textstyle \frac{1}{22} (10t_1 + 190)
\geq 0; \quad \mu_{1}(u,\chi_{2},*)  = \textstyle \frac{1}{22} (t_1+ 96) \geq 0; \\ % 2
\mu_{11}(u,\chi_{2},*) & = \textstyle \frac{1}{22} (-10t_1 + 96)
\geq 0; \quad \mu_{0}(u,\chi_{3},*)  = \textstyle \frac{1}{22} (10t_2 + 506) \geq 0; \\ % 5
\mu_{1}(u,\chi_{3},*) & = \textstyle \frac{1}{22} (t_2+ 231) \geq
0; \qquad
\mu_{11}(u,\chi_{3},*)  = \textstyle \frac{1}{22} (-10t_2 + 242) \geq 0. \\ % 8
\end{split}
\]

Case 5. Let $\chi(u^{11}) = 2\chi(2a)-\chi(2b)$. Using Proposition
\ref{P:3} for the characters $\chi_{2}$ and $\chi_{3}$ of $G$, we
obtain  the following system
\[
\begin{split}
\mu_{0}(u,\chi_{2},*) & = \textstyle \frac{1}{22} (10t_1 + 174)
\geq 0; \qquad \mu_{1}(u,\chi_{2},*)  = \textstyle \frac{1}{22} (t_1+ 112) \geq 0; \\ % 2
\mu_{11}(u,\chi_{2},*) & = \textstyle \frac{1}{22} (-10t_1 + 112)
\geq 0; \quad \mu_{0}(u,\chi_{3},*)  = \textstyle \frac{1}{22} (10t_2 + 462) \geq 0; \\ % 5
\mu_{1}(u,\chi_{3},*) & = \textstyle \frac{1}{22} (t_2+275) \geq
0; \qquad \quad
\mu_{11}(u,\chi_{3},*)  = \textstyle \frac{1}{22} (-10t_2 + 286) \geq 0. \\ % 8
\end{split}
\]

Case 6. Let $\chi(u^{11}) = -\chi(2a)+2\chi(2b)$. Using
Proposition \ref{P:3} for the characters $\chi_{2}, \chi_{3}$ and
$\chi_{4}$ of $G$, we obtain  the following system
\[
\begin{split}
\mu_{0}(u,\chi_{2},*) & = \textstyle \frac{1}{22} (10t_1 + 126)
\geq 0; \quad \mu_{11}(u,\chi_{2},*)  = \textstyle \frac{1}{22} (-10t_1+ 160) \geq 0; \\ % 2
\mu_{0}(u,\chi_{3},*) & = \textstyle \frac{1}{22} (10t_2 + 330)
\geq 0; \quad \mu_{11}(u,\chi_{3},*)  = \textstyle \frac{1}{22} (-10t_2 + 418) \geq 0; \\ % 5
\mu_{1}(u,\chi_{2},*) & = \textstyle \frac{1}{22} (t_1 + 160) \geq
0; \qquad \mu_{1}(u,\chi_{3},*)  = \textstyle \frac{1}{22} (t_2 + 407) \geq 0; \\ % 6
& \mu_{0}(u,\chi_{4},*)  = \textstyle \frac{1}{22} (120 \nu_{2a} + 200 \nu_{2b} - 10 \nu_{11a} + 798) \geq 0. \\ % 9
\end{split}
\]

Case 7. Let  $\chi(u^{11}) = -2\chi(2a)+3\chi(2b)$. Using
Proposition \ref{P:3} for the characters $\chi_{2}$ and
$\chi_{3}$ of $G$,  we obtain  the following system
\[
\begin{split}
\mu_{0}(u,\chi_{2},*) & = \textstyle \frac{1}{22} (10t_1 + 110)
\geq 0; \quad
\mu_{11}(u,\chi_{2},*)  = \textstyle \frac{1}{22} (-10t_1 + 176) \geq 0; \\ % 4
\mu_{0}(u,\chi_{3},*) & = \textstyle \frac{1}{22} (10t_2 + 286)
\geq 0; \quad
\mu_{11}(u,\chi_{3},*)  = \textstyle \frac{1}{22} (-10t_2 + 462) \geq 0; \\ % 8
\mu_{1}(u,\chi_{2},*) & = \textstyle \frac{1}{22} (t_1 + 176) \geq
0; \qquad \mu_{1}(u,\chi_{3},*)  = \textstyle \frac{1}{22} (t_2 + 451) \geq 0. \\ % 6
\end{split}
\]

Case 8. Let $\chi(u^{11}) = -3\chi(2a)+4\chi(2b)$. Using
Proposition \ref{P:3} for the characters $\chi_{2}$ and
$\chi_{3}$ of $G$, we obtain  the following system
\[
\begin{split}
\mu_{0}(u,\chi_{2},*) & = \textstyle \frac{1}{22} (10t_1 + 110)
\geq 0; \quad
\mu_{11}(u,\chi_{2},*)  = \textstyle \frac{1}{22} (-10t_1 + 192) \geq 0; \\ % 4
\mu_{0}(u,\chi_{3},*) & = \textstyle \frac{1}{22} (10t_2 + 242)
\geq 0; \quad
\mu_{11}(u,\chi_{3},*)  = \textstyle \frac{1}{22} (-10t_2 + 506) \geq 0; \\ % 8
\mu_{1}(u,\chi_{2},*) & = \textstyle \frac{1}{22} (t_1 + 192) \geq
0; \qquad \mu_{1}(u,\chi_{3},*) = \textstyle \frac{1}{22} (t_2 +
495 )\geq 0.
\end{split}
\]
In all eight cases we obtained systems of inequalities that have no solutions.

\noindent$\bullet$ Let $u$ be a unit of order $26$. By (\ref{E:1})
and Proposition \ref{P:2} we have
$$
\nu_{2a}+\nu_{2b}+\nu_{13a}+\nu_{13b}=1.
$$
Put\quad  $t_1 = 15 \nu_{2a} - \nu_{2b}$,\quad  $t_2 =
\nu_{2a}$\quad  and\quad  $t_3 = 243 \nu_{2a} + 35 \nu_{2b} - 6
\nu_{13a} + 7 \nu_{13b}$.
Since $|u^{13}|=2$ and $|u^{2}|=13$, by parts (iii) and (v) of the
Theorem we need to consider $8\cdot 18=144$ cases. We will
index them by $\chi(u^{13})$ and consider two possibilities.

First, let \qquad $\chi(u^{13})\in\{ \chi(2a),\quad 4\chi(2a)-3\chi(2b),
\quad 3\chi(2a)-2\chi(2b), \quad 2\chi(2a)-\chi(2b), \quad
-\chi(2a)+2\chi(2b),\quad -2\chi(2a)+3\chi(2b), \quad
-3\chi(2a)+4\chi(2b)\}$ \qquad and  put
$$
(\alpha,\beta)= \tiny{
\begin{cases}
(408,320),\quad & \text{if}\qquad  \chi(u^{13}) = \chi(2a) ;\\
(540,188),\quad & \text{if}\qquad  \chi(u^{13}) =4\chi(2a)-3\chi(2b) ;\\
(496,232),\quad & \text{if}\qquad  \chi(u^{13}) = 3\chi(2a)-2\chi(2b);\\
(452,276),\quad & \text{if}\qquad  \chi(u^{13}) = 2\chi(2a)-\chi(2b);\\
(320,408),\quad & \text{if}\qquad  \chi(u^{13}) = -\chi(2a)+2\chi(2b);\\
(276,452),\quad & \text{if}\qquad  \chi(u^{13}) = -2\chi(2a)+3\chi(2b);\\
(232,496),\quad & \text{if}\qquad  \chi(u^{13}) = -3\chi(2a)+4\chi(2b).\\
\end{cases}}
$$
Now using Proposition \ref{P:3} for the character  $\chi_{3}$, we
obtain  the following system
\[
\mu_{0}(u,\chi_{3},*) = \textstyle \frac{1}{26} (528t_2 +\alpha) \geq 0; \quad \mu_{13}(u,\chi_{3},*)= \textstyle \frac{1}{26} (-528t_2+\beta) \geq 0, \\ % 8
\]
which has no integral solution.

Second, suppose that   $\chi(u^{13}) = \chi(2b)$. Using
Proposition \ref{P:3} for the characters  $\chi_{2}$, $\chi_{3}$,
$\chi_{4}$ and $\chi_{31}$, we obtain the following system
\[
\begin{split}
\mu_{0}(u,\chi_{2},*)  & = \textstyle \frac{1}{26} (12t_1 + 142)\geq 0; \quad \mu_{1}(u,\chi_{2},*) = \textstyle \frac{1}{26} (-12t_1 + 144) \geq 0; \\ % 2
\mu_{0}(u,\chi_{3},*)  & = \textstyle \frac{1}{26} (528 t_2 + 364)\geq 0; \quad \mu_{13}(u,\chi_{3},*) = \textstyle \frac{1}{26} (-528 t_2 + 364) \geq 0; \\ % 8
\mu_{1}(u,\chi_{31},*) & = \textstyle \frac{1}{26} (t_3+ 93512)\geq 0; \quad \mu_{4}(u,\chi_{31},*) = \textstyle \frac{1}{26} (-t_3+ 93583) \geq 0; \\ % 100
&\mu_{13}(u,\chi_{2},*)  = \textstyle \frac{1}{26} (-180 \nu_{2a} + 12 \nu_{2b} + 144) \geq 0; \\ % 4
&\mu_{0}(u,\chi_{4},*)   = \textstyle \frac{1}{26} (144 \nu_{2a} + 240 \nu_{2b} + 800) \geq 0, \\ % 9
\end{split}
\]
which has no integral solution.

\noindent$\bullet$ Let $u$ be a unit of order $33$. By (\ref{E:1})
and Proposition \ref{P:2} we have
$$
\nu_{3a}+\nu_{3b}+\nu_{3c}+\nu_{11a}=1.
$$
Since $|u^{11}|=3$, we have to consider $104$ cases accordingly to the Appendix.
Put $t_1 = 35 \nu_{3a} + 8 \nu_{3b} -
\nu_{3c}$, $t_2 = 14 \nu_{3a} - 13 \nu_{3b} - 4 \nu_{3c} -
\nu_{11a}$ and $t_3 = 105 \nu_{3a} - 3 \nu_{3b} + 6 \nu_{3c} -
\nu_{11a}$. First, when $\chi(u^{11})$ takes values from the
first column of the following table, we have appropriate coefficients
$\alpha_i$ given in the second column
\[
\tiny
\begin{array}{|c|c|}\hline
&\\
\chi(u^{11}) & (\alpha_1, \alpha_2,\alpha_3,\alpha_4, \alpha_5, \alpha_6) \\
& \\
\hline
& \\
\chi(3a)                        & (213, 108, 346, 388, 980, 665)  \\
\chi(3b)                        & (159, 135,  361, 400, 764, 773) \\
2\chi(3a)-14\chi(3b)+13\chi(3c) & (33, 33, 47,58, 1430, 440)      \\
2\chi(3a)-13\chi(3b)+12\chi(3c) & (51,  189, 65,76, 1412, 449)    \\
2\chi(3a)-12\chi(3b)+11\chi(3c) & (69, 180, 83, 94, 1394, 458)    \\
\chi(3a)-11\chi(3b)+11\chi(3c)  & (105, 162, 496, 130, 1358, 476) \\
\chi(3a)-10\chi(3b)+10\chi(3c)  & (33, 33, 478, 166, 1160, 575)   \\
\chi(3a)-10\chi(3b)+10\chi(3c)  & (123, 153, 487, 148, 1340, 485) \\
& \\
\hline
\end{array}
\]
From the last table we obtain the system
\[
\begin{split}
\mu_{0}(u,\chi_{2},*) & = \textstyle \frac{1}{33} (20t_1+\alpha_1)
\geq 0; \qquad
\mu_{11}(u,\chi_{2},*)  = \textstyle \frac{1}{33} (-10t_1 +  \alpha_2) \geq 0; \\ % 4
\mu_{0}(u,\chi_{3},*) & = \textstyle \frac{1}{33}(-20t_2+\alpha_3)
\geq 0; \quad
\mu_{11}(u,\chi_{3},*)  = \textstyle \frac{1}{33} (10t_2+ \alpha_4) \geq 0; \\ % 8
\mu_{0}(u,\chi_{4},*) & = \textstyle \frac{1}{33} (20t_3+\alpha_5)
\geq 0; \qquad
\mu_{11}(u,\chi_{4},*)  = \textstyle \frac{1}{33} (-10t_3+ \alpha_6) \geq 0, \\ % 12
\end{split}
\]
which has no integral solution.

Finally, if $\chi(u^{11})\in\{ \chi(3c),\quad
2\chi(3a)-15\chi(3b)+14\chi(3c)\}$, then we get
\[
\begin{split}
\mu_{0}(u,\chi_{2},*) & = \textstyle \frac{1}{33} (20t_1+ 15) \geq 0; \qquad \mu_{3}(u,\chi_{2},*) = \textstyle \frac{1}{33} (-2t_1+ 15) \geq 0, \\ % 3
\end{split}
\]
and if $ \chi(u^{11}) = 2\chi(3a)-11\chi(3b)+10\chi(3c)$, we get
\[
\begin{split}
\mu_{0}(u,\chi_{2},*) & = \textstyle \frac{1}{33} (20t_1+ 87) \geq
0; \qquad
\mu_{11}(u,\chi_{2},*) = \textstyle \frac{1}{33} (-10t_1+ 171) \geq 0, \\ % 3
\end{split}
\]
both of which have no integer solutions.

\noindent$\bullet$ Let $u$ be a unit of order $35$. By (\ref{E:1})
and Proposition \ref{P:2} we have
$$
 \nu_{5a}+\nu_{5b}+\nu_{7a}=1.
$$
Put $t_1 = 8 \nu_{5a} + 3 \nu_{5b} + 3 \nu_{7a}$ and $t_2 =
\nu_{5a} - 4 \nu_{5b}$. Since $|u^{7}|=5$, by part (vi) of the
Theorem we have to consider nine cases.

Case 1. Let $\chi(u^{7}) = \chi(5a)$, then
\[
\begin{split}
\mu_{0}(u,\chi_{2},*) & = \textstyle \frac{1}{35} (24 t_1 + 193)
\geq 0; \qquad \mu_{7}(u,\chi_{2},*)  = \textstyle \frac{1}{35} (-6t_1 + 153) \geq 0; \\ % 4
\mu_{0}(u,\chi_{3},*) & = \textstyle \frac{1}{35} (-24 t_2+ 360)
\geq 0; \quad \mu_{7}(u,\chi_{3},*)  = \textstyle \frac{1}{35} (6t_2+ 365) \geq 0; \\ % 8
&\qquad  \mu_{7}(u,\chi_{2},2) = \textstyle \frac{1}{35} (30 \nu_{5a} + 12 \nu_{7a} + 103) \geq 0. \\ % 124
\end{split}
\]

Case 2. Let $\chi(u^{7}) = \chi(5b)$, then
\[
\begin{split}
\mu_{0}(u,\chi_{2},*) & = \textstyle \frac{1}{35} (24 t_1+ 173)
\geq 0; \qquad
\mu_{7}(u,\chi_{2},*)  = \textstyle \frac{1}{35} (-6t_1+ 158) \geq 0; \\ % 4
\mu_{0}(u,\chi_{3},*) & = \textstyle \frac{1}{35} (-24t_2+ 380)
\geq 0; \quad
\mu_{7}(u,\chi_{3},*) = \textstyle \frac{1}{35} (6 t_2+ 360) \geq 0; \\ % 8
&\qquad  \mu_{0}(u,\chi_{4},*)  = \textstyle \frac{1}{35} (240 \nu_{5a} + 72 \nu_{7a} + 798) \geq 0.\\ % 9
\end{split}
\]

Case 3. Let $\chi(u^{7}) = 5\chi(5a)-4\chi(5b)$, then
\[
\begin{split}
\mu_{0}(u,\chi_{2},*) & = \textstyle \frac{1}{35} (24t_1+ 273)
\geq 0; \qquad
\mu_{7}(u,\chi_{2},*)  = \textstyle \frac{1}{35} (-6t_1+ 133) \geq 0; \\ % 4
\mu_{0}(u,\chi_{3},*) & = \textstyle \frac{1}{35} (-24t_2+ 280)
\geq 0; \quad
\mu_{7}(u,\chi_{3},*) = \textstyle \frac{1}{35} (6 t_2+ 385) \geq 0. \\ % 8
\end{split}
\]

Case 4. Let $\chi(u^{7}) = 4\chi(5a)-3\chi(5b)$, then
\[
\begin{split}
\mu_{0}(u,\chi_{2},*) & = \textstyle \frac{1}{35} (24t_1+ 253)
\geq 0; \qquad
\mu_{7}(u,\chi_{2},*)  = \textstyle \frac{1}{35} (-6t_1+ 138) \geq 0; \\ % 4
\mu_{0}(u,\chi_{3},*) & = \textstyle \frac{1}{35} (-24 t_2+ 300)
\geq 0; \quad
\mu_{7}(u,\chi_{3},*)  = \textstyle \frac{1}{35} (6t_2+ 380) \geq 0; \\ % 8
&\qquad \mu_{0}(u,\chi_{2},3)  = \textstyle \frac{1}{35} (96 \nu_{5a} - 24 \nu_{5b} + 24 \nu_{7a} + 146) \geq 0. \\ % 153
\end{split}
\]

Case 5. Let $\chi(u^{7}) = 3\chi(5a)-2\chi(5b)$, then
\[
\begin{split}
\mu_{0}(u,\chi_{2},*) & = \textstyle \frac{1}{35} (24t_1+ 233)
\geq 0; \qquad
\mu_{7}(u,\chi_{2},*) = \textstyle \frac{1}{35} (-6t_1+ 143) \geq 0; \\ % 4
\mu_{0}(u,\chi_{3},*) & = \textstyle \frac{1}{35} (-24t_2+ 320)
\geq 0; \quad
\mu_{7}(u,\chi_{3},*)  = \textstyle \frac{1}{35} (6t_2+ 375) \geq 0; \\ % 8
&\mu_{0}(u,\chi_{2},3)  = \textstyle \frac{1}{35} (96 \nu_{5a} - 24 \nu_{5b} + 24 \nu_{7a} + 126) \geq 0. \\ % 153
\end{split}
\]

Case 6. Let $\chi(u^{7}) = 2\chi(5a)-\chi(5b)$, then
\[
\begin{split}
\mu_{0}(u,\chi_{2},*) & = \textstyle \frac{1}{35} (24t_1+ 213)
\geq 0; \qquad
\mu_{7}(u,\chi_{2},*)  = \textstyle \frac{1}{35} (-6t_1+ 148) \geq 0; \\ % 4
\mu_{0}(u,\chi_{3},*) & = \textstyle \frac{1}{35} (-24t_2+ 340)
\geq 0; \quad
\mu_{7}(u,\chi_{3},*)  = \textstyle \frac{1}{35} (6t_2+ 370) \geq 0; \\ % 8
&\mu_{0}(u,\chi_{7},2)  = \textstyle \frac{1}{35} (72 \nu_{5a} - 48 \nu_{5b} + 24 \nu_{7a} + 676) \geq 0. \\ % 133
\end{split}
\]

Case 7. Let $\chi(u^{7}) = -\chi(5a)+2\chi(5b)$, then
\[
\begin{split}
\mu_{0}(u,\chi_{2},*) & = \textstyle \frac{1}{35} (24t_1+ 153)
\geq 0; \qquad \mu_{7}(u,\chi_{2},*)  = \textstyle \frac{1}{35}
(-6t_1+ 163) \geq 0; \\  \mu_{0}(u,\chi_{3},*) & = \textstyle
\frac{1}{35} (-24t_2+ 400) \geq 0; \quad
\mu_{7}(u,\chi_{3},*)  = \textstyle \frac{1}{35} (6 t_2+ 355) \geq 0. \\ % 8
\end{split}
\]

Case 8. Let $\chi(u^{7}) = -2\chi(5a)+3\chi(5b)$, then
\[
\begin{split}
\mu_{0}(u,\chi_{2},*) & = \textstyle \frac{1}{35} (24t_1+ 133)
\geq 0; \qquad
\mu_{5}(u,\chi_{2},*)  = \textstyle \frac{1}{35} (-6t_1+ 112) \geq 0; \\ % 3
\mu_{0}(u,\chi_{3},*) & = \textstyle \frac{1}{35} (-24t_2+ 420)
\geq 0; \quad
\mu_{7}(u,\chi_{3},*)  = \textstyle \frac{1}{35} (6t_2+ 350) \geq 0; \\ % 8
& \mu_{0}(u,\chi_{2},3)  = \textstyle \frac{1}{35} (96 \nu_{5a} - 24 \nu_{5b} + 24 \nu_{7a} + 26) \geq 0. \\ % 153
\end{split}
\]

Case 9. Let $\chi(u^{7}) = -3\chi(5a)+4\chi(5b)$, then
\[
\begin{split}
\mu_{0}(u,\chi_{2},*) & = \textstyle \frac{1}{35} (24t_1+ 113)
\geq 0; \qquad
\mu_{5}(u,\chi_{2},*)  = \textstyle \frac{1}{35} (-6t_1+ 92) \geq 0; \\ % 3
\mu_{0}(u,\chi_{3},*) & = \textstyle \frac{1}{35} (-24t_2+ 440)
\geq 0; \quad
\mu_{7}(u,\chi_{3},*)  = \textstyle \frac{1}{35} (6t_2+ 345) \geq 0; \\ % 8
&\quad  \mu_{0}(u,\chi_{2},3)  = \textstyle \frac{1}{35} (96 \nu_{5a} - 24 \nu_{5b} + 24 \nu_{7a} + 6) \geq 0. \\ % 153
\end{split}
\]
In all of the above cases we obtained systems that have no integer solutions.

\noindent$\bullet$ Let $u$ be a unit of order $39$. By (\ref{E:1})
and Proposition \ref{P:2} we have
$$
\nu_{3a}+\nu_{3b}+\nu_{3c}+\nu_{13a}+\nu_{13b}=1.
$$
Put $t_1 = 35 \nu_{3a} + 8 \nu_{3b} -  \nu_{3c}$, $t_2 = 14
\nu_{3a} - 13 \nu_{3b} - 4 \nu_{3c}$, $t_3 = 35 \nu_{3a} -
\nu_{3b} + 2 \nu_{3c}$ and $t_4 =9 \nu_{3c} - 7 \nu_{13a} + 6
\nu_{13b}$.
Since $|u^{13}|=3$ and $|u^{3}|=13$, by part (vi) of the Theorem and the
Appendix we have to consider $1872$ cases. Using the LAGUNA
package \cite{LAGUNA}, we find out that in all cases we have
the system of inequalities which has no integral solutions
\[
\begin{split}
\mu_{0}(u,\chi_{2},*) & = \textstyle \frac{1}{39} (24t_1+ 213)
\geq 0; \qquad
\mu_{1}(u,\chi_{2},*) = \textstyle \frac{1}{39} (t_1+108) \geq 0; \\ % 2
\mu_{13}(u,\chi_{2},*) & = \textstyle \frac{1}{39} (-12
t_1+108)\geq 0; \quad
\mu_{0}(u,\chi_{3},*) = \textstyle \frac{1}{39} (-24 t_2+ 336) \geq 0; \\ % 5
\mu_{1}(u,\chi_{3},*) & = \textstyle \frac{1}{39} (-t_2 + 378)\geq
0; \qquad
\mu_{13}(u,\chi_{3},*)  = \textstyle \frac{1}{39} (12 t_2+ 378) \geq 0; \\ % 8
\mu_{0}(u,\chi_{4},*) & = \textstyle \frac{1}{39} (72t_3 + 990)
\geq 0; \qquad \mu_{1}(u,\chi_{31},*)  = \textstyle \frac{1}{39}
(-t_4 + 93548) \geq 0; \\
&\qquad \mu_{3}(u,\chi_{31},*)  = \textstyle \frac{1}{39} (2t_4 + 93548) \geq 0. \\ % 96
\end{split}
\]

\noindent$\bullet$ Let $u$ be a unit of order $55$. By (\ref{E:1})
and Proposition \ref{P:2} we have
$$
\nu_{5a}+\nu_{5b}+\nu_{11a}=1.
$$
Put $t_1 = 8 \nu_{5a} + 3 \nu_{5b}$\quad and\quad  $t_2 = \nu_{5a}
- 4 \nu_{5b} -  \nu_{11a}$. Since $|u^{11}|=5$, by part (iv) of
the Theorem we have to consider nine cases that we collect into
four groups.

Case 1. \quad When \;
$\chi(u^{11}) \in \{ \; \chi(5a), \; \chi(5b), \;
4\chi(5a)-3\chi(5b), \; 3\chi(5a)-2\chi(5b), \\
2\chi(5a)-\chi(5b), \; -\chi(5a)+2\chi(5b)\; \}$
then we put
\[
\begin{split}
(\alpha,\beta)= \tiny{
\begin{cases}
(175,135) & \text{if }\quad  \chi(u^{11}) = \chi(5a);\\%1
(155,140) & \text{if }\quad  \chi(u^{11}) = \chi(5b);\\%2
(235,120) & \text{if } \quad \chi(u^{11}) = 4\chi(5a)-3\chi(5b);\\%4
(215,125) & \text{if }\quad  \chi(u^{11}) = 3\chi(5a)-2\chi(5b);\\%5
(195,130) & \text{if }\quad  \chi(u^{11}) = 2\chi(5a)-\chi(5b);\\%6
(135,145) & \text{if }\quad  \chi(u^{11}) = -\chi(5a)+2\chi(5b).\\%7
\end{cases}}
\end{split}
\]
Using Proposition \ref{P:3} for the character
$\chi_{2}$ of $G$, we get the system
\[
\begin{split}
\mu_{0}(u,\chi_{2},*) & = \textstyle \frac{1}{55} (40t_1 + \alpha)
\geq 0; \quad
\mu_{11}(u,\chi_{2},*)  = \textstyle \frac{1}{55} (-10t_1+\beta) \geq 0; \\ % 4
& \qquad \mu_{1}(u,\chi_{2},*)  = \textstyle \frac{1}{55} (t_1+\beta) \geq 0, \\ % 2
\end{split}
\]
which has no integral solution.

In the remaining three cases below we have no integral solutions as well.

Case 2. Let\quad  $\chi(u^{11}) = 5\chi(5a)-4\chi(5b)$, then
\[
\begin{split}
\mu_{0}(u,\chi_{2},*) & = \textstyle \frac{1}{55} (40 t_1 + 255)
\geq 0; \quad
\mu_{11}(u,\chi_{2},*)  = \textstyle \frac{1}{55} (-10 t_1 + 115) \geq 0; \\ % 4
\mu_{1}(u,\chi_{2},*) & = \textstyle \frac{1}{55} (t_1+ 115) \geq
0; \qquad
\mu_{0}(u,\chi_{3},*)  = \textstyle \frac{1}{55} (-40t_2+ 290) \geq 0; \\ % 5
\mu_{11}(u,\chi_{3},*) & = \textstyle \frac{1}{55} (10 t_2+ 395)
\geq 0; \quad
\mu_{1}(u,\chi_{3},*)  = \textstyle \frac{1}{55} (-t_2+ 384) \geq 0. \\ % 6
\end{split}
\]

Case 3. Let\quad  $\chi(u^{11}) = -2\chi(5a)+3\chi(5b)$, then
\[
\begin{split}
\mu_{0}(u,\chi_{2},*) & = \textstyle \frac{1}{55} (40t_1+ 115)
\geq 0; \qquad
\mu_{11}(u,\chi_{2},*)  = \textstyle \frac{1}{55} (-10t_1+ 150) \geq 0; \\ % 4
\mu_{1}(u,\chi_{2},*) & = \textstyle \frac{1}{55} (t_1+ 150) \geq
0; \qquad\quad
\mu_{0}(u,\chi_{3},*)  = \textstyle \frac{1}{55} (-t_2+ 430) \geq 0; \\ % 5
\mu_{1}(u,\chi_{3},*) & = \textstyle \frac{1}{55} (-t_2+ 349) \geq
0; \qquad
\mu_{11}(u,\chi_{3},*)  = \textstyle \frac{1}{55} (10 t_2+ 360) \geq 0. \\ % 8
\end{split}
\]

Case 4. Let $\chi(u^{11}) = -3\chi(5a)+4\chi(5b)$, then
\[
\begin{split}
\mu_{0}(u,\chi_{2},*) & = \textstyle \frac{1}{55} (40t_1+ 95) \geq
0; \qquad
\mu_{11}(u,\chi_{2},*) = \textstyle \frac{1}{55} (-10t_1+ 155) \geq 0; \\ % 4
\mu_{1}(u,\chi_{2},*) & = \textstyle \frac{1}{55} (t_1+ 155) \geq
0; \qquad\quad
\mu_{11}(u,\chi_{3},*)  = \textstyle \frac{1}{55} (10t_2+ 355) \geq 0; \\ % 8
\mu_{0}(u,\chi_{3},*) & = \textstyle \frac{1}{55} (-40t_2+ 450)
\geq 0; \quad
\mu_{1}(u,\chi_{3},*)  = \textstyle \frac{1}{55} (-t_2+ 344) \geq 0. \\ % 6
\end{split}
\]

\noindent$\bullet$ Let $u$ be a unit of order $65$. By (\ref{E:1})
and Proposition \ref{P:2} we have
$$
    \nu_{5a}+\nu_{5b}+\nu_{13a}+\nu_{13b}=1.
$$
Put $t_1 = 8 \nu_{5a} + 3 \nu_{5b}$, $t_2 =  \nu_{5a} - 4
\nu_{5b}$ and $t_3 = 10 \nu_{5a} - 6 \nu_{13a} + 7 \nu_{13b}$.
Since $|u^{13}|=5$ and $|u^{5}|=13$, by parts  (iv) and (v)  of
the Theorem we have to consider  $162$ cases.\\
First, let $(\alpha,\beta)= \tiny
\begin{cases}
(175, 135) &\text{if}\quad \chi(u^{13}) = \chi(5a);\\
(155,40)&\text{if}\quad \chi(u^{13}) = \chi(5b);\\
(255,115)&\text{if}\quad  \chi(u^{13}) = 5\chi(5a)-4\chi(5b);\\
(135,145)&\text{if}\quad  \chi(u^{13}) = -\chi(5a)+2\chi(5b);\\
(115,150)&\text{if}\quad  \chi(u^{13}) = -2\chi(5a)+3\chi(5b);\\
(95,155)&\text{if}\quad   \chi(u^{13}) = -3\chi(5a)+4\chi(5b).\\
\end{cases}
$
\\
Then the following pair of inequalities has no integral solution
\[
\begin{split}
\mu_{0}(u,\chi_{2},*) & = \textstyle \frac{1}{65} (48t_1+\alpha)
\geq 0; \qquad
\mu_{13}(u,\chi_{2},*)  = \textstyle \frac{1}{65} (-12t_1+\beta) \geq 0. \\
\end{split}
\]
In the remaining cases we put
\[
(\alpha_1, \alpha_2, \alpha_3,\alpha_4)= \tiny
\begin{cases}
(235,120,380,300),\quad & \text{if }\quad \chi(u^{13}) =
4\chi(5a)-3\chi(5b);\\
(215,125,375,320),\quad & \text{if }\quad \chi(u^{13}) =
3\chi(5a)-2\chi(5b);\\
(195, 130, 370, 340),\quad & \text{if }\quad \chi(u^{13}) =
2\chi(5a)-\chi(5b),\\
\end{cases}
\]
and also we parametrize $(\alpha_5,\alpha_6)$ by values of
$\chi(u^{13})$ and $\chi(u^{5})$ accordingly to the following table
\[
\tiny
\begin{array}{|c|c|c|c|c|}\hline
& & & &\\
& \chi(u^{13})=& \chi(u^{13})=& \chi(u^{13}) = &  \chi(u^{5})=\\
& 4\chi(5a)-3\chi(5b) & 3\chi(5a)-2\chi(5b) & 2\chi(5a)-\chi(5b) &\\
& & & &\\
\hline
& & & &\\
&(93508, 93708) & (93518, 93668) &  (93528, 93628) &  \chi(13a)\\
&(93521, 93721) & (93531, 93681) &  (93541, 93641) &  \chi(13b)\\
&(93404, 93604) & (93414, 93564) &  (93424, 93524) &  9\chi(13a)-8\chi(13b)\\
&(93417, 93617) & (93427, 93577) &  (93437, 93537) &  8\chi(13a)-7\chi(13b)\\
&(93430, 93630) & (93440, 93590) &  (93450, 93550) &  7\chi(13a)-6\chi(13b)\\
&(93443, 93643) & (93453, 93603) &  (93463, 93563) &  6\chi(13a)-5\chi(13b)\\
&(93456, 93656) & (93466, 93616) &  (93476, 93576) &  5\chi(13a)-4\chi(13b)\\
&(93469, 93669) & (93479, 93629) &  (93489, 93589) &  4\chi(13a)-3\chi(13b)\\
(\alpha_5, \alpha_6)   &(93482, 93682) & (93492, 93642) &  (93502, 93602) &  3\chi(13a)-2\chi(13b)\\
&(93495, 93695) & (93505, 93655) &  (93515, 93615) &  2\chi(13a)-\chi(13b)\\
&(93534, 93734) & (93544, 93694) &  (93554, 93654) &  -\chi(13a)+2\chi(13b)\\
&(93547, 93747) & (93557, 93707) &  (93567, 93667) &  -2\chi(13a)+3\chi(13b)\\
&(93560, 93760) & (93570, 93720) &  (93580, 93680) &  -3\chi(13a)+4\chi(13b)\\
&(93573, 93773) & (93583, 93733) &  (93593, 93693) &  -4\chi(13a)+5\chi(13b)\\
&(93586, 93786) & (93596, 93746) &  (93606, 93706) &  -5\chi(13a)+6\chi(13b)\\
&(93599, 93799) & (93609, 93759) &  (93619, 93719) &  -6\chi(13a)+7\chi(13b)\\
&(93612, 93812) & (93622, 93772) &  (93632, 93732) &  -7\chi(13a)+8\chi(13b)\\
&(93625, 93825) & (93635, 93785) &  (93645, 93745) &  -8\chi(13a)+9\chi(13b)\\
& & & &\\
\hline
\end{array}
\]
In all of these cases we obtain the following system
\[
\begin{split}
\mu_{0}(u,\chi_{2},*) & = \textstyle \frac{1}{65} (48
t_1+\alpha_1) \geq 0; \qquad
\mu_{13}(u,\chi_{2},*)  = \textstyle \frac{1}{65} (-12 t_1+\alpha_2) \geq 0; \\ %4
\mu_{0}(u,\chi_{3},*) & = \textstyle \frac{1}{65} (-48 t_2+
\alpha_3) \geq 0; \quad
\mu_{13}(u,\chi_{3},*)  = \textstyle \frac{1}{65} (12t_2+ \alpha_4) \geq 0; \\ % 8
\mu_{1}(u,\chi_{31},*) & = \textstyle \frac{1}{65} (t_3+\alpha_5)
\geq 0; \qquad \mu_{10}(u,\chi_{31},*) = \textstyle \frac{1}{65}
(-4t_3+\alpha_6) \geq 0, \\ % 73
\end{split}
\]
that has no integral solutions.

\noindent$\bullet$ Let $u$ be a unit of order $77$. By (\ref{E:1})
and Proposition \ref{P:2} we have $\nu_{7a}+\nu_{11a}=1$.
Then we obtain the system
\[
\begin{split}
\mu_{0}(u,\chi_{2},*) & = \textstyle \frac{1}{77} (180 \nu_{7a} +
161) \geq 0; \quad
\mu_{0}(u,\chi_{2},2)  = \textstyle \frac{1}{77} (-120 \nu_{7a} + 98) \geq 0, \\ % 81
\end{split}
\]
which has no integral solution.

\noindent$\bullet$ Let $u$ be a unit of order $91$. By (\ref{E:1})
and Proposition \ref{P:2} we have
$$
\nu_{7a}+\nu_{13a}+\nu_{13b}=1.
$$
Since $|u^{7}|=13$, by part   (v)  of the Theorem we have to
consider $18$ cases. But in all cases we obtain the same non-compatible
system of inequalities
\[
\begin{split}
\mu_{0}(u,\chi_{2},*) & = \textstyle \frac{1}{91} (216 \nu_{7a} +
161) \geq 0; \qquad
\mu_{13}(u,\chi_{2},*)  = \textstyle \frac{1}{91} (-36 \nu_{7a} + 140) \geq 0. \\ % 4
\end{split}
\]

\noindent$\bullet$ Let $u$ be a unit of order $143$. By
(\ref{E:1}) and Proposition \ref{P:2} we have
$$
\nu_{11a}+\nu_{13a}+\nu_{13b}=1.
$$
Since $|u^{11}|=13$, by part   (v)  of the Theorem we have to
consider $18$ cases. But in all cases we obtain the same system of
inequalities
\[
\begin{split}
\mu_{0}(u,\chi_{3},3)&=\textstyle\frac{1}{143}(120\nu_{11a}+88)\geq
0; \qquad
\mu_{0}(u,\chi_{4},*)  = \textstyle \frac{1}{143} (-120 \nu_{11a} + 770) \geq 0, \\ % 5
\end{split}
\]
which has no integral solutions.

\vspace{15pt}
\centerline{{\bf Appendix.}}
\vspace{10pt}
\centerline{Possible partial augmentations $( \nu_{3a}, \nu_{3b}, \nu_{3c} )$ for units of order 3:}
\[
\begin{split}
\; \{ \; ( -3, \nu_{3b}, \nu_{3c} ) \; & | \quad \; 5 \le \nu_{3b} \le 7,\quad \nu_{3a}+\nu_{3b}+\nu_{3c}=1 \; \} \; \cup \\
\; \{ \; ( -2, \nu_{3b}, \nu_{3c} ) \; & | \quad \; 1 \le \nu_{3b} \le 11,\quad \nu_{3a}+\nu_{3b}+\nu_{3c}=1 \; \} \; \cup \\
\; \{ \; ( -1, \nu_{3b}, \nu_{3c} ) \; & | \; -3 \le \nu_{3b} \le 14,\quad \nu_{3a}+\nu_{3b}+\nu_{3c}=1 \; \} \; \cup \\
\; \{ \; ( \;\; 0, \nu_{3b}, \nu_{3c} ) \; & | \; -7 \le \nu_{3b} \le 16,\quad \nu_{3a}+\nu_{3b}+\nu_{3c}=1 \; \} \; \cup \\
\; \{ \; ( \;\; 1, \nu_{3b}, \nu_{3c} ) \; & | \; -11 \le \nu_{3b} \le 12,\quad \nu_{3a}+\nu_{3b}+\nu_{3c}=1 \; \} \; \cup \\
\; \{ \; ( \;\; 2, \nu_{3b}, \nu_{3c} ) \; & | \; -15 \le \nu_{3b} \le 8,\quad \nu_{3a}+\nu_{3b}+\nu_{3c}=1 \; \}. \\
\end{split}
\]

\bibliographystyle{plain}
\bibliography{Tex_Suz_Kimm}

\end{document}